\documentclass[11pt,letter,times]{article}

\usepackage{times}
\usepackage{fullpage}
\usepackage{latexsym}
\usepackage{amsmath,amsthm,amsfonts,amssymb,amstext}

\usepackage{ifpdf}
\ifpdf
\usepackage[pdftex]{graphicx}
\else
\usepackage[dvips]{graphicx}
\fi

\ifpdf
\usepackage[pdftex,colorlinks,linkcolor=blue,citecolor=blue]{hyperref}
\else
\fi

\def\draft{0}   

\ifnum\draft=1 
    \def\ShowAuthNotes{1}
\else
    \def\ShowAuthNotes{0}
\fi

\newcommand{\newitem}{\vspace{-1.5ex}\item}

\newcommand{\Ex}[2]{\underset{#1}{\mathbb{E}} [#2] }

\newtheorem{theorem}{Theorem}
\newtheorem{lemma}[theorem]{Lemma}

\theoremstyle{definition}
\newtheorem{definition}[theorem]{Definition}

\newcommand{\etal}{{\it et al }}
\newcommand{\e}{\varepsilon}

\ifnum\ShowAuthNotes=1
  \newcommand{\authnote}[2]{{ \bf [#1's Note: #2]}}
\else
  \newcommand{\authnote}[2]{}
\fi

\newcommand{\half}{\frac{1}{2}}

\title{Routing Complexity of Faulty Networks
   \ifnum\draft=1 \\
         { \small {\sc Working Draft }}
   \fi }

\author{Omer Angel\footnote{Orsay University, Paris. \texttt{omer.angel@math.u-psud.fr}} \footnote{Research done while at the Weizmann Institute.}
 \and Itai Benjamini\footnote{ The Weizmann Institute of Science.  Rehovot 76100
 Israel.
\texttt{\{itai.benjamini,eran.ofek,udi.wieder\}@weizmann.ac.il}}
\and Eran Ofek\footnotemark[\value{footnote}] \and Udi
Wieder\footnotemark[\value{footnote}]}

\begin{document}

\begin{titlepage}
\maketitle
\thispagestyle{empty}

\begin{abstract}
One of the fundamental problems in distributed computing is how to
efficiently perform routing in a faulty network in which each link
fails with some probability. This paper investigates how big the
failure probability can be, before the capability to efficiently
find a path in the network is lost. Our main results show tight
upper and lower bounds for the failure probability which permits
routing, both for the hypercube and for the $d-$dimensional mesh.
We use tools from percolation theory to show that in the
$d-$dimensional mesh, once a giant component appears
--- efficient routing is possible. A different behavior is
observed when the hypercube is considered. In the hypercube there
is a range of failure probabilities in which short paths exist
with high probability, yet finding them must involve querying
essentially the entire network. Thus the routing complexity of the
hypercube shows an asymptotic phase transition. The critical
probability with respect to routing complexity lies in a different
location then that of the critical probability with respect to
connectivity. Finally we show that an oracle access to links (as
opposed to local routing) may reduce significantly the complexity
of the routing problem. We demonstrate this fact by providing
tight upper and lower bounds for the complexity of routing in the
random graph $G_{n,p}$.
\end{abstract}

%
%
\end{titlepage}
%
\newpage

\section{Introduction}
The goal of this paper is to investigate the effectiveness of
routing in faulty networks. Suppose that a network is represented
by a graph $G$. Two kinds of fault models are common in the
theoretical literature: Worst case faults and random faults which
are our concern. In the random fault model it is assumed that each
component of the network fails with some probability and
independently of all other components. In this paper we consider
edge failures so we assume each edge in $G$ fails independently
with some probability $q=1-p$. Let $G_p$ be the resulting graph.
One can ask what is the probability that two nodes $u$ and $v$
remain connected in $G_p$. This had been the focus of much
research concerning the existence of {\em giant components} in
such graphs, and the critical values of $p$ for the existence of
those, cf. \cite{aks82,ns86,karlin94,abs03,ksv02}. But in many
applications the fact that a path between u and v exists is not
sufficient, one wants to be able to {\em find} the path in a
distributed manner.



It is known that if the topology of a graph has some randomness,
then the existence of short paths in a graph does not guarantee
the ability of efficiently finding them. For instance a cycle with
a random matching has a logarithmic diameter \cite{bolobas88}, yet
paths connecting a given pair of nodes can not be found in less
than $\sqrt n$ time \cite{smallworld}. This phenomenon is
especially acute when considering `natural' networks such as the
world wide web, social networks, P2P networks etc, in which
typically the network size is huge, the diameter of the network is
small and the challenge is to find short paths within a time
complexity that is comparable to the diameter. Indeed Kleinberg's
model of the small world phenomenon \cite{smallworld,kleinberg01},
is aimed at explaining the ability to \emph{find} short paths in
social networks (and not merely their existence). In the context
of P2P several randomized topologies were proposed along with
routing algorithms that find short paths in the random graph cf.
\cite{skipgraph,gummadi03,mnw04}. In the context of P2P networks,
many routing algorithm are able to find paths between nodes even
when nodes or links fail cf. \cite{hk03,nw03,chord}. While our
findings do not apply directly to these networks, we expect that
are main result would hold for them as well. See Section
\ref{sec:results} for more details.

\medskip

In this paper we analyze the algorithmic complexity of finding a
path between nodes $u,v$ in $G_p$ as a function of the failure
probability. In particular we seek to find the exact values of $p$
for which it is possible to perform routing in $G_p$ within time
complexity that is comparable to the diameter. One difficulty is
that with positive probability $u,v$ are in distinct components of
$G_p$. We therefore restrict our attention to the case where a
giant component exists, and condition on the event that $u,v$ are
connected.

Our findings present a complex picture. We show that for some
graphs, as the $d-$dimensional mesh, efficient routing is always
possible, i.e. it is easy to find with high probability short
paths between nodes within the giant component (whenever it
exists). However, for other graphs, such as the hypercube,
efficient routing is possible only for some failure probabilities.
In other words, there is a range of failure probabilities for
which with high probability a giant component exists, the diameter
of the giant component is small, yet in order to find a path
between nodes it is necessary to probe a large portion of the
graph. We provide tight upper and lower bounds on the routing
complexity, indicating the exact location of the transition.


\subsection{The Model}

\begin{definition}
Given a graph $G_p$ and two vertices $u,v$, a \emph{routing
algorithm} is an algorithm that is allowed to \emph{probe} whether
an edge exists in $G_p$, and outputs a path between $u,v$ if such
exists. A routing algorithm is said to be \emph{local}, if the
first edge it probes is adjacent to $u$ and subsequently it probes
only edges to (an end point of) which it has already established a
path from $u$.
\end{definition}

Local algorithms aim to capture the realistic constraints of
routing in a network. If each node is a server in a network and
$u$ wishes to send a message to $v$ then $u$ must find a path to
$v$ while probing edges it has already reached. In Section
\ref{sec: OracleLocalRouting} we show that a local router may
require \emph{exponentially} more probes than a non local one,
thus the distinction between the two kinds of algorithms is
necessary. A non local routing algorithm may be referred to as an
\emph{oracle} routing algorithm.
Denote by $\{u\sim v\}$ the event that $u$ is indeed connected to
$v$.

\begin{definition}
Given a graph $G$, probability $p$ and a routing algorithm $A$,
the \emph{routing complexity} of $A$ denoted by $comp(A)$, with
respect to the nodes $u,v$, is the random variable that counts the
queries $A$ makes (i.e.\ edges probed) to find a path between
$u,v$ in $G_p$, conditioned on $\{u\sim v\}$.
\end{definition}

The routing complexity measures how many \emph{probes} are needed
to route a message from $u$ to $v$ in $G_p$, assuming this routing
is possible. We do not consider here any computations that the
algorithm performs. As indicated above, the question is most
interesting when $\Pr[u\sim v]$ is bounded away from zero, and
indeed we limit our discussion to this case. A simple upperbound
on the routing complexity could be achieved by performing a BFS
search on $G_p$. In terms of the routing complexity this is
tantamount to probing the entire graph. However there may exist
algorithms which achieve a much smaller routing complexity. In
particular, if the diameter of $G$ is small, we are interested in
finding a routing algorithm with a complexity that is comparable
to the \emph{actual distance} between the nodes, or show that none
exists.

We stress that the routing complexity measures the expected
complexity of finding a path between two specified vertices, and
does not necessarily indicate the difficulty of performing a full
blown routing scheme. Small routing complexity may be seen as the
minimal requirement of fault tolerance in networks. Naturally such
a weak requirement strengthens our lower bounds and weakens our
upper bounds. In this paper we focus on analyzing the hypercube
and the $d-$dimensional grid, which are probably the most widely
investigated topologies in this context.
\subsection{Related Work}
Denote by $H_{n,p}$ the $n-$dimensional hypercube, when each edge
is deleted with probability $1-p$ and survives with probability
$p$. Random subgraphs of the hypercube had been the focus of much
research. It is known (see eg.\ \cite{es}) that if $p<\half$ then
with high probability $H_{n,p}$ is not connected and if $p>\half$
then with high probability $H_{n,p}$ is connected. A classic
result by Ajtai, Komlos and Szemeredi \cite{aks82} states that if
$p\ge n^{-1}(1+\e)$ for any fixed $\e>0$ then with high
probability\footnote{Throughout this paper, the term `with high
probability' means with probability that tends to $1$ as
$n\to\infty$.} $H_{n,p}$ contains a giant component (i.e. a
component with $\Theta(2^n)$ nodes), while if $\e < 0$ then w.h.p
a giant component will not exist. This result was sharpened by
Bollob\'{a}s \etal in \cite{bkl92} and then by Borgs \etal in
\cite{bshss03}.

A related notion to routing complexity is that of
\emph{emulation}. Roughly speaking, network $A$ emulates network
$B$ if $A$ can perform \emph{any} computation $B$ performs with a
constant slowdown. When the emulating network is a random subgraph
the notion of emulation implies not only that short paths could be
found but also that they do not create bottlenecks in the
computation. Hastad \etal \cite{hastad87,hastad89} considered
\emph{node} failures, and showed that if $p$ is a constant close
enough to $1$, then $H_{n,p}$ could \emph{emulate} $H_n$ with a
small slowdown. Cole \etal \cite{cole95} proved that a faulty
butterfly network can perform efficient permutation routing even
if each node or edge fails with some constant probability.
Emulation under worst case faults were considered by Leighton
\etal \cite{leighton98}. In particular these results imply that if
the failure probability is small enough then it is possible to
find paths efficiently between nodes in the giant component. On
the other hand Angel and Benjamini \cite{ab03} showed that if
$p<\frac{1}{\sqrt n}$ then the hypercube could not be embedded in
its giant component with constant distortion. This result suggests
that when $\frac{1}{n}<p<\frac{1}{\sqrt n}$ then even though a
giant component exists w.h.p, it defines a metric that is
different then that of the hypercube.

\medskip

Let $M^d_p$ be a $d$-dimensional mesh with $M^d$ nodes, in which
each edge is deleted with probability $1-p$. It is known that for
each $d$ there exists a critical probability $p^d_c$ such that if
$p<p_c^d$ then w.h.p.\ there will not be a giant component in
$M^d_p$, and if $p>p_c^d$ then w.h.p.\ $M^d_p$ will contain a
giant component. The exact values of the critical probabilities
are not always known. It is known that $p_c^2=\half$ and that
$p_c^d=(1+o(1))/2d$ and is decreasing in $d$. See the book by
Grimmett \cite{percolation} and the references therein.
\authnote{Udi}{Find an estimate on the size of the g.c. as a
function of d,p} Kaklamanis \etal \cite{kklmrrt} showed that if
$p$ is large enough then $M^2_p$ can emulate $M^d$ with $O(\log
n)$ slowdown. Mathies \cite{mathies} extended this result for any
$p>p^2_c=\half$. These results do not imply an efficient routing
algorithm. Naor and Wieder \cite{naor03scalable} used planar
duality to prove that efficient routing is possible in $M^2_p$
whenever $p>\half$. Cole \etal \cite{cole93} proved that a two
dimensional array can tolerate a constant fraction of worst case
faults and still emulate the non faulty array with a constant
slowdown.

\subsection{Summary of Results}\label{sec:results}

Recall that $H_{n,p}$ denotes a random subgraph of the
$n-$dimensional hypercube obtained by selecting each edge
independently with probability $p$. As mentioned it is known that
when $p>(1+\e)n^{-1}$ with high probability a giant component
exists (\cite{aks82}). Furthermore it is implicit in the proof,
that the diameter of the component is polynomial in $n$. If
however $p\le n^{-1}(1-\e)$ then the size of the largest connected
component is $o(2^n)$ w.h.p. This suggests that efficient routing
in the giant component of $H_{n,p}$ might be possible for any
$p\ge n^{-1}(1+\e)$. Our work shows that this is \emph{not} the
case: there is a threshold for efficient routing that lies in a
\emph{different} location. The following theorem provides a
complete characterization of the routing complexity as a function
of the failure probability:

\begin{theorem}\label{thm:hc}
For a fixed $\alpha$, let $p=n^{-\alpha}$.
\begin{description}
\newitem $(i)$
If $\alpha > 1/2+\beta$ for $\beta>0$, any local routing algorithm
in the hypercube $H_{n,p}$, makes at least $2^{\Omega(n^\beta)}$
queries w.h.p..
\newitem $(ii)$
There is a local routing algorithm $A$ on the hypercube  such that
the following holds: For any $\alpha<1/2$ there exists
$k=k(\alpha)$ so that for any two vertices, $comp(A)<n^k$ with
probability $1-\exp(-cn^{1-\alpha})$ for some constant $c>0$.

\end{description}
\end{theorem}

Thus for $p=n^{-\alpha}$ for $\half<\alpha<1$, an intriguing
phenomenon occurs: the giant component of $H_{n,p}$ shares some
structural properties of $H_n$, in particular it has diameter
$poly(n)$ (w.h.p.), and has roughly the same expansion of $H_n$,
yet the ability to find short paths is lost. Angel and Benjamini
proved in \cite{ab03} that for these failure probabilities $H_n$
could not be embedded in $H_{n,p}$ with constant distortion, so
the result of Part $(i)$ is not entirely surprising, yet the
techniques we use are different then that of \cite{ab03}. Part
$(i)$ is proven by showing that if $p$ is small enough then a ball
centered at $v$ is likely to look more or less like a tree rooted
at $v$, which contains closed edges. Now, in order to reach $v$
from $u$ it is necessary to find a leaf which is connected to $v$
via an open path, an event which is proven to be rare. Our
technique is general enough to be used on other families of
graphs.

It is proven in \cite{ab03} that if $\alpha<1/2$ then there is an
embedding of $H_{n}$ in $H_{n,p}$ with constant distortion. This
embedding is used to derive the matching upperbound of part
$(ii)$. Note that the algorithm of $(ii)$ does not depend on
$\alpha$, and only its efficiency changes. Therefore if $\alpha =
0$, i.e. there are no faults, then the algorithm reduces to a
greedy algorithm which routes along the hypercube's shortest
paths.

Many popular P2P topologies share some structural similarities
with the hypercube cf. \cite{chord,skipgraph,pastry}. We did not
prove that Theorem \ref{thm:hc} holds for these topologies, yet it
is reasonable to assume that that this is the case. If so then the
Theorem implies that if the network suffers many faults, flooding
and gossiping techniques would remain efficient means to locate
data (in terms of latency) while the routing based exact search
algorithms fail.

\medskip

The phenomenon described in Theorem \ref{thm:hc} does not occur in
all graphs. Recall that $M^d_p$ denotes a $d$-dimensional mesh
with $M^d$ nodes, in which each edge is deleted with probability
$1-p$.

\begin{theorem}\label{thm:mesh}
Let $u,v$ be two vertices at distance $n$ in $M^d$. There exists a
local routing algorithm in $M^d_p$ so that if $p>p_c^d$ then the
expected routing complexity is $O(n)$.
\end{theorem}

Thus in the mesh when $p$ is large enough so that a giant
component appears it is possible to find paths between two
vertices in time comparable to the distance between the nodes. It
is important to note that if we allow $p$ to be close enough to
$1$ then Theorem \ref{thm:mesh} is fairly easy prove. The main
difficulty is proving that the Theorem's statement is correct for
\emph{any $p>p_c^d$}, this involves some deep results from
Percolation Theory.

\medskip

The previous two theorems assumed the routing algorithms are
local, i.e. they are only allowed to probe edges for which they
have already established a path. What if we remove the locality
assumption and allow the routing algorithm to probe \emph{any}
edge in the graph, we call this model \emph{oracle} routing. On
first glance it might seem as if oracle routing may not change
considerably the routing complexity. Yet, in Section \ref{sec:
OracleLocalRouting} we show a graph in which there is an
exponential gap between the routing complexity with respect to
oracle routing and that of local routing. We also provide tight
upper and lower bounds for routing in $G_{n,p}$, and show that in
this natural model oracle routing outperforms local routing.

\medskip

In the next section a lemma which provides a lower bound for
routing complexity in a general scenario is proved. In
Sections~\ref{sec:hypercube} and \ref{sec:mesh} we prove our
results for the hypercube and the mesh respectively. Section
\ref{sec: OracleLocalRouting} concerns the oracle routing model.
Finally Section \ref{sec:conjecture} discusses some related open
problems.

\section{The Lower Bound Lemma} \label{sec: lower_bound_lemma}
In this Section we prove a lemma which is instrumental in proving
hardness of local routing on various graphs. The basic intuition
could be seen through the following example: Consider a graph in
which there are exactly $d$ edge disjoint paths of length $2$
between nodes $u,v$. Now assume each edge remains open with
probability $\frac{1}{\sqrt d}$. We expect that both $u$ and $v$
would be connected to about $\sqrt{d}$ open edges, thus by the
birthday paradox w.h.p there would be an open path of length $2$
between the nodes. Assuming such a path exists, it is easy to see
that $\Omega(d)$ edges should be probed w.h.p before one of these
paths is found.

%

This intuition could be generalized as follows:  If $S$ is a
subset of the nodes, and $v\in S$ while $u \not \in S$, then a
path from $u$ to $v$ must at some point find an edge in the cut
$(S,\bar S)$ which is connected to $v$. If the probability that an
edge in the cut is connected to $v$ via the set $S$ is low enough,
then many such edges should be probed before a path is found. More
formally: for a set $S$ and vertices $u,v\in S$ we write $\{(x\sim
y)\in S\}$ for the event that\ $x$ is connected to $y$ by an open
path in the set $S$. Similarly $\{(x\sim e)\in S\}$ denotes the
event that $x$ is connected to an end point of the edge $e$ via a
path in the set $S$.

\begin{lemma} \label{lem: localhard}
Let $V=S\cup\bar S$ be a disjoint partition of the vertex set of a
graph and $v\in S$ a vertex. Assume for any edge $e$ crossing the
cut $(S,\bar S)$ we have $\Pr[(v\sim e)\in S] \le\eta$, and let
$X$ be the number of queries made by a local routing algorithm
from $u$ to $v$, then
\begin{equation}\label{eq:localhard}
\Pr[X<t] \le \frac{t\eta + \Pr[(u\sim v) \in S]}{\Pr[(u\sim v)]}.
\end{equation}
If $u\in\bar S$ then the numerator becomes $t\eta$.
\end{lemma}
%

\authnote{Omer}{It is possible to add the expectation formulation
as a corollary with the last bit of the original proof.}

\begin{proof}
If $((u \sim v) \in S)^c$ (which is always the case if $u\in\bar
S$), then any path from $u$ to $v$ crosses the cut. Finding a path
from $u$ to $v$ involves finding a path in $S$ from some edge of
the cut to $v$. Each probed edge of the cut has probability at
most $\eta$ of being in such a path. For any set of $t$ edges in
the cut $(S,\bar S)$, the probability that one of them is
connected to $v$ in $S$ is at most $t\eta$. Thus the probability
of finding a path from $u$ to $v$ by probing only those $t$ edges
of the cut is at most $t\eta+\Pr[(u \sim v) \in S]$.

Since this bound is uniform in the set of edges, the fact that
edges may be chosen adaptively does not invalidate the bound.
Similarly, the constraint of only probing edges reachable from $u$
only reduces the possible sets of queries. Finally, the
denominator stems from conditioning on $\{u\sim v\}$ (in the
complexity definition).
\end{proof}

\subsection{An Illustrative Example}
In following we use Lemma \ref{lem: localhard} to lowerbound the
routing complexity of the double binary tree. The double binary
tree of depth $n$ denoted $TT_n$ is constructed by taking two
binary trees of uniform depth $n$ and identifying their
corresponding leaves. Let $x,y$ be the two roots of the trees.
First we identify the
failure probabilities for which $x,y$ are connected with
probability which is bounded away from $0$.
\
\begin{lemma}
If $\frac{1}{\sqrt 2} < p \le 1$ then there exists a path between
$x$ and $y$ in $TT_{n,p}$ with probability bounded away from 0. If
$p\le \frac{1}{\sqrt 2}$ then w.h.p.\ no such paths exists.
\end{lemma}

\begin{proof}
In order for a path to exist between $x$ and $y$, it must be the
case that there exists an open branch from a leaf $w$ to the root
of the first tree $x$, and that the mirroring branch from $w$ to
$y$ is also open. This is equivalent to the case of a single tree
where each edge is open with probability $p^2$. It is well known
that the critical probability of a Galton Watson tree (or the
binary branching process) is $\half$. See for instance
\cite{stirzaker} for details.
\end{proof}


Now we show that for \emph{any} $p<1$, the local routing
complexity is exponential in the diameter of the graph.
\begin{theorem}\label{thm:doubletree}
Let $\frac{1}{\sqrt 2} < p < 1$. For some $c>0$ and any $a,n$, any
local router between the two roots of $TT_n$ makes at least
$ap^{-n}$ queries with probability at least $1-ca$.
\end{theorem}

\authnote{Omer}{rephrase this to hold for any graph with
neighborhood of $y$ that is a tree?}

\begin{proof}
Apply Lemma~\ref{lem: localhard} with $S$ being the second tree to
get the desired bound: Clearly we may have $\eta=p^n$. The nodes
$x$ and $y$ can be connected only via the cut $(S,\bar{S})$, this
happens with probability at least $c(p)$. Lemma \ref{lem:
localhard} now implies
\[\Pr[A < ap^{-n}] < \frac{ap^{-n} p^{n}}{c(p)} = \frac{a}{c(p)}\]
\end{proof}
If we set $a$ to be a decaying function (say $\frac{1}{ n}$) then
the probability a local router would probe less than
$\frac{p^{-n}}{n}$ is $O(\frac{1}{n})$.
The double binary tree has the interesting property that an
\emph{oracle} routing algorithm may find a path between $x$ and
$y$ with a \emph{polynomial} number of probes. See Section
\ref{sec: OracleLocalRouting}.

\section{Hypercube --- Tight Upper and Lower Bounds} \label{sec:hypercube}
In this section we show the exact location of the probability $p$,
in which the routing complexity shows a phase transition between
being exponential and being polynomial (in $n$). The idea is to
show that when $p<\frac{1}{\sqrt n}$ then balls around nodes look
more-or-less like trees, and therefore when trying to reach node
$v$, a routing algorithm would need to `penetrate' a tree through
its leaves, as was demonstrated in the double binary tree graph.
When $p>\frac{1}{\sqrt n}$ then there are enough edges so that
some variant of greedy routing will find a path within polynomial
time.
\subsection{The Lower Bound}

Here we apply Lemma \ref{lem: localhard} to get a
lower bound on the local routing complexity of the hypercube when
$p<n^{-1/2}$. 
The given bound translates to a fractional exponential (in $n$)
bound on the routing complexity.

\begin{proof}[Proof of Theorem~\ref{thm:hc}(i)]
We apply Lemma~\ref{lem: localhard} to the hypercube with $S$
being a ball of radius $l=n^\beta$ around $y$, for some
$0<\beta<\alpha-1/2$.

The first stage is to bound $\eta$ of the lemma, i.e. bound the
probability $v$ is connected to an edge on the boundary of $S$ via
a path within $S$ . We show that for large $n$, for any $e$
connecting $S$ and $\bar S$ we have $\Pr[(v\sim e)\in S] \le \eta$
holds with $\eta = 2n^{(\beta-\alpha)n^{\beta}}$.
Let $x$ be the endpoint of $e$ in $S$ with $d(x,v)=l$. Consider a
path from $v$ to $x$ in $S$ as a sequence of coordinates in which
consecutive steps are taken. Let $A_k$ be the set of such paths of
length $l+2k$ (by parity this catches all paths).

For $k=0$ we have $|A_0| = l!$ since a path of $A_0$ uses each of
the $l$ coordinates exactly once. To bound $|A_k|$
we show a map from $A_k$ to $A_{k-1}$ that maps at most $n\cdot
l^2$ paths to each path. Existence of such a map implies
$|A_k|\leq nl^2|A_{k-1}|$ and therefore by induction $|A_k| \le
n^k l^{2k} l!$. To define the map, consider the first $l+1$ steps
of a path. Since the path remains in the ball $S$, at least one of
the coordinates is repeated. Take such a repeated coordinate and
eliminate its first two occurrences. It is easy to see that this
maps a path $p\in A_k$ to a path $p'\in A_{k-1}$. To reconstruct
$p$ from $p'$ one needs to know which coordinate was removed ($n$
possibilities) and the indices at which it appeared
($\binom{l+1}{2} \le l^2$ possibilities). Thus the pre-image of
$p'$ contains at most $n l^2$ paths from $A_k$.

This bound clearly counts many paths more than once, as well as
many non-simple paths, but it is good enough. Each simple
path in $A_k$ is open with probability $p^{l+2k}$, and so
\begin{align*}
\Pr[(v\sim x) \in S]
  &\le  \sum_{k=0}^{\infty} p^{l+2k} n^k l^{2k} l! \\
  &\le (lp)^l \sum_{k=0}^{\infty} (nl^2p^2)^k
   = \frac{n^{(\beta-\alpha)n^\beta}}{1-n^{2\beta+1-2\alpha}} .
\end{align*}
For large $n$ the denominator is close to 1, hence
$\eta=2n^{(\beta-\alpha)n^\beta}$ is a valid choice.

Next, we estimate the other terms in \eqref{eq:localhard}. Since
each of $u$ and $v$ is in the giant component with probability
tending to 1, $\Pr[(u \sim v)] \to 1$. If $u\not\in S$ then
$\Pr[(u\sim v)\in S]=0$. Otherwise, suppose $d(u,v)=m\le l$. The
same argument as above shows that the number of paths in $S$ of
length $m+2k$ from $u$ to $v$ is at most $m! (nl^2)^k$ and hence
the probability that any of them are open satisfies
\[
\Pr[(u\sim v)\in S] \le \sum_{k=0}^{\infty} p^{m+2k} n^k l^{2k} m!
                     =  \frac{m! p^m}{1-n^{2\beta+1-2\alpha}} .
\]
The denominator tends to 1 and for $m\le l$, the numerator is
$o(1)$ because $mp \le lp =n^{\beta-\alpha}$.

Using Lemma~\ref{lem: localhard}, we now see that if the
complexity of a local router in the hypercube is $A$, then
\[
\Pr[A< n^{(\alpha-\beta)n^{\beta}}/n]
  \le \frac{2/n+ \Pr[(u\sim v) \in S]}{\Pr[(u\sim v)]}
  \to 0.
\]
\end{proof}

\subsection{The Upper Bound}
Next we show that when $p$ is large, local routing on the
hypercube may be performed using $n^k$ probes with high
probability. This is a variation on the result of \cite{ab03}
showing that in this regime the metric distortion of the
percolation is bounded. This shows that there is indeed an
asymptotic phase transition in the complexity of routing on the
hypercube. The proof below shows that $k=O((1-2\alpha)^{-1})$,
though it would be interesting to know the exact dependence of $k$
on $\alpha$ (the optimal $k$ need not be integral).

\begin{proof}[Proof of Theorem~\ref{thm:hc}(ii)]
Here, the terms neighbor and distance relate to the metric of the
hypercube before percolation. Percolation neighbor and percolation
distance are used for the percolated hypercube $H_{n,p}$.\\

We refer to the definition of a good vertex from \cite{ab03},
which roughly means having a high degree in $H_{n,p}$. The
condition that a vertex is good is determined by the neighborhood
of percolation radius $2$ around it. In \cite{ab03}, Section $(2)$
the following is proved:
\begin{description}
\item (1) Any given vertex is good with probability  $1-\exp(-cn^{1-\alpha})$.
\item (2) With probability $1-\exp(-cn)$, all pairs of good vertices at distance up to 3 have
percolation distance at most $l$ for some
$l=l(\alpha)=O((1-2\alpha)^{-1})$.
\end{description}

Now the algorithm is straight forward. Pick arbitrarily a path
from $u$ to $v$, of minimal length: $u=u_0,u_1,\ldots,u_m=v$, and
use BFS iteratively to find a path from $u_i$ to $u_{i+1}$. With
probability tending to 1 all the vertices of the path, including
$u$ and $v$ are good (each one is not good with probability at
most $\exp(-cn^{1-\alpha})$, there are at most $n$ vertices in the
path). On this event, the percolation distance between $u_i$ and
$u_{i+1}$ is at most $l$ and a path from $u_i$ to $u_{i+1}$ can be
found by, say, BFS of complexity $n^l$. The total complexity is at
most $n^{l+1}$.
\end{proof}
\paragraph{Remark:} A natural approach would be to use greedy routing,
i.e. at each routing step, probe edges that reduce the Hamming
distance to the target. While this strategy may work most of the
way, in the final steps a more extensive search is required. It
may be the case though that a greedy approach at the early stages
of the routing would reduce the exponent in the complexity of the
algorithm.


\section{The Mesh --- Upper Bound} \label{sec:mesh}
In this section we show that the phenomenon observed for
hypercubes does not apply when the mesh is considered, i.e.
whenever a giant component exists, it is possible to efficiently
route between nodes. Consider a cube of the $d$-dimensional mesh,
i.e. a submesh with $M^d$ nodes, and let each edge remain open
with probability with some fixed $p$, and be closed with
probability $q=1-p$. Let $d(\cdot,\cdot)$ denote distance in the
mesh, and $D(\cdot,\cdot)$ denote the distance in the giant
component (which may be referred as percolation distance). We seek
a path between two vertices $u,v$ in the cube with $d(u,v)=n$ (the
cube size is $M^d$ which may be much larger than $n$). We are
interested in the routing complexity in terms of $n$ when $p$ is
fixed. As mentioned, there exists a number $p_c^d$ such that if
$p\le p^d_c$ then $\Pr[u\sim v] = o(1)$ as $n\to\infty$, so
hereinafter we assume $p>p^d_c$. For such $p$ there is a giant
cluster in the cube, and with probability bounded from 0, both $u$
and $v$ are in the giant component and therefore connected.

We give an algorithm that efficiently finds a short path from $u$
to $v$. The case of $d=2$ was solved by Naor and Wieder in
\cite{naor03scalable}, where planar duality is used to show that
in a two dimensional grid with $n^2$ vertices, the routing
complexity is $O(n)$ w.h.p. It is important to note that it is
fairly easy to find a path between $u,v$ if we assume that $p$ is
sufficiently close to $1$. The main difficulty is pushing the
probability $p$ all the way down to $p_c^d$. In order to do that
we need some fairly recent and strong results from Percolation
Theory.

\subsection{The Routing Algorithm}

The idea of the algorithm is as follows. Consider $n$ vertices
which belong to some shortest path between $u,v$. With high
probability many of them are in the giant component and the
percolation distance between them is not too large. The algorithm
searches around each of them, until the next one is found. More
formally:
\begin{enumerate}
\newitem Fix $u=u_0,u_1,u_2,\ldots,u_n=v$ to be a shortest path
between $u$ and $v$. Start from $u_0=u$.
\newitem Assume $u_i$ has been reached. Exhaustively probe edges
around $u_i$ (using say BFS) until some vertex $u_j$ with $j>i$ is
reached.
\newitem Repeat at most $n$ times until reaching $u_n=v$.
\end{enumerate}

Note that a BFS up do distance $k$ from a vertex takes only
$O(k^d)$ queries since only edges of the mesh at distance $k$ from
the starting point may be reached. The key point is that it is
very unlikely at any iteration that a large depth is needed.
Correctness of the algorithm is clear since the search at each
stage stops once a closer approximation to $v$ is found. If an
open path from $u$ to $v$ exists, then some path will be found.


\begin{proof}[Proof of Theorem~\ref{thm:mesh}]
Let $u_i$ be some vertex along the chosen path to $v$ that is in
the giant component. Let $u_j$ be the next vertex along the path
in the giant component. It follows that the $j-i-1$ vertices along
the path between them are outside the giant component, an event
that is exponentially unlikely (see \cite{percolation}):
\[
\Pr[|j-i| > k] < e^{-c_1 k} \qquad \text{for some $c_1=c_1(p)>0$.}
\]

\authnote{Omer}{It would be nice if we had an exact reference that
in a finite cube this is true for any set of $k$ vertices (or at
least for a monotone path of length $k$.)}

Note that $j$ always exists since $v$ is assumed to be in the
giant component. In practice, $u_j$ might be skipped over by the
algorithm if some further vertex $u_k$ is reached first. If the
algorithm explores a neighborhood of $u_i$, it finds a further
vertex of the path at distance at most $D(u_i,u_j)$. Thus to bound
the number of queries the algorithm makes to reach some $u_k$ we
use the following Lemma, which is a proper restatement of result
by Antal and Pisztora \cite{ap96, gm90}.

\begin{lemma} \label{lem:new_ap}
For any $p>p^d_c$ and any $x,y$ in a cube $M^d$ of the infinite
mesh, let $D(x,y)$ be the percolation distance (in $M^d$) between
them. For some $\rho,c_2>0$ depending only on the dimension and
$p$, and for any $a>\rho\cdot d(x,y))$
\[
\Pr[(D(x,y) > a) \wedge (x\sim y)\in M^d] < e^{-c_2 a} .
\]
\end{lemma}

\authnote{Omer}{We should fix this, at least in the journal
version. either we add a note on how to get away from the
boundary, or we show that the boundary does not invalidate A-P.
The alternative would be to use regular A-P for
vertices at least polylog from the boundary, and have another
stage at the beginning to get from $u$ away from the boundary (and
another at the end). When making a large step the argument used
for the slabs shows that the boundary is not significant.}

Either $d(u_u,u_j)$ is large or it is small. In the latter case,
$D(u_i,u_j)$ is unlikely to be large, and the former case is
itself unlikely:
\begin{align*}
\Pr[D(u_i,u_j) > k]
 &< \Pr[d(u_i,u_j) > k/\rho] +
    \Pr[(d(u_i,u_j) \le k/\rho) \wedge(D(u_i,u_j)>k)] \\
 &< e^{-c_1 k/\rho} + e^{-c_2 k} < e^{-c_3k} .
\end{align*}
Consequently, if $A_i$ is the number of queries made from $u_i$,
\[
\Pr[A_i > k] < \Pr[D(u_i,u_j) > ck^{1/d}] < e^{-c_4k^{1/d}} .
\]
Since this is summable, for each vertex $u_i$ of the path that is
in the giant component the expected work to get from $u_i$ to a
further vertex is $O(1)$.

The number of queries made by the algorithm is at most the sum
over all vertices of the path in the giant component of the work
to progress from them (actually it is less since some may be
skipped over, and some queries may be duplicated). By additivity
of expectation, $\Ex{}{\sum A_i}=O(n)$.
\end{proof}

\section{Oracle Routing vs. Local Routing}\label{sec: OracleLocalRouting}

In this section we consider routing algorithms that are allowed to
query \emph{any} edge, and not just edges to which it has
established a path. This is called \emph{oracle routing}.
Surprisingly, it might be the case that a huge gap exists between
the complexity of local and oracle routing. A simple (yet somewhat
artificial) example for this is the double binary tree $TT_n$ with
fixed $\sqrt{1/2}<p<1$. In section \ref{sec: lower_bound_lemma} we
showed that any local routing algorithm which finds a path between
the two roots of $TT_n$ w.h.p.\ makes exponentially many queries.
The following theorem shows that oracle routing algorithms can do
significantly better.

\begin{theorem}
There is an oracle router between the two roots of $TT_n$ with
average complexity $cn$ for some $c=c(p)<\infty$ and any
$p>\sqrt{1/2}$.
\end{theorem}

\begin{proof}
A simple path between the two roots is just a branch up to level
$n$ in the first tree joined to the corresponding branch in the
second tree. The oracle router is very simple: To find a path
from the root to level $n$ that is open in both trees, query edges
together with their corresponding edges in the second tree. Each
such pair of edges is open with probability $p^2 > 1/2$. The
problem is equivalent to finding a path from the root to level $n$
in a super-critical Galton-Watson tree, and a depth first search
accomplishes this in expected complexity linear in $n$.

To see this, observe that any branch of the infinite binary
branching process (the Galton-Watson tree) that fails to reach
level $n$ has expected size $c(p)$ (and is in fact exponentially
unlikely to be large), see \cite{lp}. Since at most $n$ bad
branches are encountered before reaching level $n$, the routing
complexity is bounded by a sum of $n$ random variables with finite
expectation and second moments, and hence is linear in $n$.
\end{proof}

A more natural example is the graph $G_{n,p}$: For each pair of
nodes $u,v$ the edge $(u,v)$ is open with probability $p$. In our
setting it could be thought of as a faulty complete graph. It
turns out that local routers can not do much better than querying
all the edges:

\begin{theorem}
Any local routing algorithm for the $G_{n,p}$ model where $p=c/n$
(for $c>1$) has an expected local routing complexity of at least
$\Omega(n^2)$.
\end{theorem}

\begin{proof}
Assume we wish to route from $u$ to $v$, and let $X$ be the number
of queries required. Let $U_t$ be the set of vertices of the graph
which are connected to $u$ by paths known to be open after $t$
queries. Thus $U_0=\{u\}$. Each vertex in $U_t$ has probability
$p$ of being connected to $v$, thus the probability of finding a
route while $|U_t|\le k$ is at most $pk$.

To reach an additional vertex given that a set of vertices has
been reached, the only option is to probe an edge connecting $U_t$
to its complement. By symmetry all such edges are equivalent, and
each has probability $c/n$ of connecting to a new vertex. Thus
$|U_t|-1$ is just a sum of $0-1$ random variables with expectation
$p=c/n$.

Since $u,v$ are both in the giant component with probability at
least some $a>0$, it follows that
\begin{align*}
\Pr[X<k]
 &< \frac{\Pr[U_k>n\sqrt{k}p]+\Pr[X<k|U_k\leq n\sqrt{k}p]}{\Pr[u\sim v]}\\
 &< \frac{\Pr[U_k>n\sqrt{k}p] + p\cdot n\sqrt{k}p}{\Pr[u\sim v]}\\
 &< \frac{\sqrt{k}/n + \frac{c^2\sqrt{k}}{n}}{a} = O(\sqrt{k}/n).
\end{align*}
where the last inequality follows from Markov's inequality. This
is close to 0 for $k=o(n^2)$ and shows that the average complexity
is $\Omega(n^2)$.
\end{proof}

Next we give tight upper and lower bounds on the oracle routing
complexity. The next theorem implies that oracle routing in this
case is better than local routing by a factor of exactly $\sqrt
n$.

\begin{theorem}
There exists a routing algorithm with average complexity
$O(n\sqrt{n})$; Any algorithm succeeds with $an\sqrt{n}$ queries
with probability at most $O(ca^{2/3} + cn^{-1})$.
\end{theorem}

\begin{proof}
Let $U_t$ and $V_t$ be the sets of vertices reachable from $u$ and
$v$ after $t$ queries. For the upper bound, consider the following
algorithm
\begin{description}
\newitem (1) Whenever there are unqueried edges between $V_t$ and
$U_t$, probe one of them,
\newitem (2) Otherwise, pick the smaller of $U_t,V_t$ and probe an
unprobed edge connecting it to a previously unreached vertex.
\newitem (3) If no such edge exists, return that $u \not\sim v$.
\end{description}

The algorithm is trivially correct. Since $U_t$ and $V_t$ grow by
one vertex at a time, they are roughly of equal size. Since each
edge is open with probability $c/n$, on average a connection
between $U_t$ and $V_t$ will be found when
$|U_t|=|V_t|=\theta(\sqrt{n})$. Since adding a vertex to either of
the sets requires a number of queries with geometric distribution
and mean $n/c$, it takes $O(n^{3/2})$ queries to find a path from
$u$ to $v$.

For the lower bound, note that $|U_t \cup V_t| \le 2+s_t$ where
$s_t$ is the number of open edges found by time $t$, and
$\Ex{}{s_t}=ct/n$. Before a connection from $u$ to $v$ is found
there must be an open (unprobed) edge between $U_t$ and $V_t$, and
the probability of that is at most $\frac{c|U_v|\cdot|V_t|}{n} \le
\frac{c(s_t+2)^2}{4n}$. Thus for any algorithm $A$ and any
$\lambda$
\[
\Pr[comp(A)<t] \le \Pr[s_t > \lambda] + \frac{c(\lambda+2)^2}{4n}
         \le \frac{ct}{n\lambda} + \frac{2c(\lambda^2+4)}{4n} .
\]
If $t=an^{3/2}$, then setting $\lambda=t^{1/3}$ results in
\[
\Pr[comp(A)<an^{3/2}] \le \frac{3 ca^{2/3}}{2} + \frac{2}{n} .
\]
\end{proof}

\section{Open Questions} \label{sec:conjecture}

So far we observed that sometimes efficient routing is possible
whenever the giant component exists, and sometimes the routing
complexity has a phase transition at a different value of the
percolation parameter. It is natural to assume that this
phenomenon relates to the growth rate of the graph. In particular:
\begin{itemize}
\item Prove or refute: there exists a family of constant degree graphs in which: the diameter
is logarithmic in the number of nodes, and the locations of the
phase transition of percolation and routing coincide (at a
location bounded away from $1$).
\end{itemize}
In particular it would be interesting to analyze De-Bruijn graphs,
Shuffle-Exchange graphs, Butterflies and other families often used
in the context of parallel computing.

It would be interesting to see hardness results for oracle
routers. Above we see that in the complete graph the best oracle
router has complexity $\theta(n^{3/2})$ where the giant component
has diameter $\theta(\log n)$. If $p$ were a small power of $n$,
then the diameter would be $O(1)$ and the complexity would still
be some power of $n$, and thus there is no bound for the
complexity in terms of the diameter. The results of \cite{ab03}
suggest that oracle routing would \emph{not} help in the
hypercube.
\begin{itemize}
\item Prove that for $\frac{1}{n}<p<\frac{1}{\sqrt n}$ the \emph{oracle} routing complexity of the
hypercube is exponential in $n$.
\end{itemize}


\begin{thebibliography}{16}

\bibitem{aks82}
M.~Ajtai, J.~Komlos, and E.~Szemeredi.
\newblock Largest random component of a k-cube.
\newblock {\em Combinatorica}, 2(1):1--7, 1982.

\bibitem{abs03}
N.~Alon, I.~Benjamini, and A.~Stacey.
\newblock Percolation on finite graphs and isoperimetric inequalities.
\newblock In {\em Annals of Probability to appear}, 2003.

\bibitem{ab03}
O.~Angel and I.~Benjamini.
\newblock A phase transition for the metric distortion of percolation on the
  hypercube.
\newblock {\em arXiv:math.PR/0306355}.

\bibitem{ap96}
P.~Antal and A.~Pisztora.
\newblock On the chemical distance in supercritical bernoulli percolation.
\newblock {\em The Annals of Probability}, (24):1036--1048, 1996.

\bibitem{skipgraph}
J.~Aspnes and G.~Shah.
\newblock Skip graphs.
\newblock In {\em Proc. 14th ACM-SIAM Symp. on Discrete Algorithms (SODA
  2003)}, pages 384--393, Jan. 2003.

\bibitem{bolobas88}
B.~Bollobas and F.~Chung.
\newblock The diameter of a cycle plus a random matching.
\newblock {\em {SIAM} Journal on Discrete Mathematics}, 1:328--333, 1988.

\bibitem{bkl92}
B.~Bollob\'{a}s, Y.~Kohayakawa, and T.~Luczak.
\newblock The evaluation of random subgraphs of the cube.
\newblock {\em Random Structures and Algorithms}, 1(3):55--90., 1992.

\bibitem{bshss03}
C.~Borgs, J.~T. Chayes, R.~van~der Hofstad, G.~Slade, and
J.~Spencer.
\newblock Random subgraphs of finite graphs: Iii. the phase transition for the
  $n$-cube.
\newblock In {\em ArXive Article math.PR/0401071}.

\bibitem{cole93}
R.~Cole, B.~M. Maggs, and R.~K. Sitaraman.
\newblock Multi-scale self-simulation: a technique for reconfiguring arrays
  with faults.
\newblock In {\em {ACM} Symposium on Theory of Computer Science (STOC)}, pages
  561--572, 1993.

\bibitem{cole95}
R.~Cole, B.~M. Maggs, and R.~K. Sitaraman.
\newblock Routing on butterfly networks with random faults.
\newblock In {\em {IEEE} Symposium on Foundations of Computer Science}, pages
  558--570, 1995.

\bibitem{es}
P.~Erdos and J.~Spencer.
\newblock Evolution of the $n$-cube.
\newblock {\em Comput. Math. Appl.}, (5):33--39., 1979.

\bibitem{percolation}
G.~Grimmett.
\newblock {\em Percolation}.
\newblock Springer-Verlag, second edition, 1999.

\bibitem{gm90}
G.~Grimmett and J.~M. Marstrand.
\newblock The supercritical phase of percolation is well behaved.
\newblock {\em Proceedings of the Royal Society (London)}, A(430):439--457,
  1990.

\bibitem{stirzaker}
G.~R. Grimmett and D.~R. Stirzaker.
\newblock {\em Probability and Random Processes}.
\newblock Oxford Science Publications, second edition, 1993.

\bibitem{gummadi03}
K.~P. Gummadi, R.~Gummadi, S.~D. Gribble, S.~Ratnasamy,
S.~Shenker, and
  I.~Stoica.
\newblock The impact of {DHT} routing geometry on resilience and proximity.
\newblock In {\em Proc. ACM SIGCOMM 2003}, pages 381--394, 2003.

\bibitem{hastad87}
J.~Hastad, T.~Leighton, and M.~Newman.
\newblock Reconfiguring a hypercube in the presence of faults.
\newblock In {\em Proceedings of the Nineteenth Annual {ACM} Symposium on
  Theory of Computing (STOC)}, pages 274--284, May 1987.

\bibitem{hastad89}
J.~Hastad, T.~Leighton, and M.~Newman.
\newblock Fast computation using faulty hypercubes.
\newblock In {\em Proceedings of the 21st Annual {ACM} Symposium on Theory of
  Computing (STOC)}, pages 251--263, 1989.

\bibitem{hk03}
K.~Hildrum and J.~Kubiatowicz.
\newblock Asymptotically efficient approaches to fault-tolerance in
  peer-to-peer networks.
\newblock In {\em 17th International Symposium on Distributed Computing
  (DISC)}, pages 321--336, 2003.

\bibitem{kklmrrt}
C.~Kaklamanis, A.~R. Karlin, F.~T. Leighton, V.~Milenkovic,
P.~Raghavan,
  S.~Rao, C.~D. Thomborson, and A.~Tsantilas.
\newblock Asymptotically tight bounds for computing with faulty arrays of
  processors.
\newblock In {\em The 31st Annual Symposium on Foundations of Computer Science
  (FOCS)}, pages 285--296, 1990.

\bibitem{karlin94}
A.~R. Karlin, G.~Nelson, and H.~Tamaki.
\newblock On the fault tolerance of the butterfly.
\newblock In {\em Proceedings of the 26th Annual ACM Symposium on the Theory of
  Computing (STOC)}, pages 125--133, 1994.

\bibitem{smallworld}
J.~Kleinberg.
\newblock The {S}mall-{W}orld phenomenon: An algorithmic perspective.
\newblock In {\em Proceedings of the 32nd ACM Symposium on Theory of
  Computing}, 2000.

\bibitem{kleinberg01}
J.~Kleinberg.
\newblock The small-world phenomenon: An algorithmic perspective.
\newblock {\em Advances in Neural Information Processing Systems (NIPS)}, (14),
  2001.

\bibitem{ksv02}
M.~Krivelevich, B.~Sudakov, and V.~Vu.
\newblock Sharp threshold for network reliability.
\newblock {\em Combinatorics, Probability and Computing}, (11):465--474, 2002.

\bibitem{leighton98}
F.~T. Leighton, B.~M. Maggs, and R.~K. Sitaraman.
\newblock On the fault tolerance of some popular bounded-degree networks.
\newblock {\em SIAM Journal on Computing}, 27(5):1303--1333., 1998.

\bibitem{lp}
R.~Lyons and Y.~Peres.
\newblock {\em Probability on Trees and Networks}.

\bibitem{mnw04}
G.~S. Manku, M.~Naor, and U.~Wieder.
\newblock Know thy neighbor's neighbor: The power of lookahead in randomized
  p2p networks.
\newblock In {\em Proceedings of the 36th ACM Symposium on Theory of Computing
  (STOC)}, 2004.

\bibitem{mathies}
T.~R. Mathies.
\newblock Percolation theory and computing with faulty arrays of processors.
\newblock In {\em Proceedings of the Third Annual Symposium on Discrete
  Algorithms (SODA)}, pages 100--103, 1992.

\bibitem{naor03scalable}
M.~Naor and U.~Wieder.
\newblock Scalable and dynamic quorum systems.
\newblock In {\em {ACM} Conf. on Principles of Distributed Computing (PODC)},
  2003.

\bibitem{nw03}
M.~Naor and U.~Wieder.
\newblock A simple fault tolerant distributed hash table.
\newblock In {\em Second International Workshop on P2P Systems (IPTPS)}, 2003.

\bibitem{ns86}
C.~M. Newman and L.~S. Schulman.
\newblock One dimensional $1/|j-i|^s$ percolation models: The existence of a
  transition for $s\leq 2$.
\newblock {\em Communications in Mathematical Physics}, 180:483--504, 1986.

\bibitem{pastry}
A.~Rowstron and P.~Druschel.
\newblock Pastry: Scalable, decentralized object location, and routing for
  large-scale peer-to-peer systems.
\newblock {\em Lecture Notes in Computer Science}, 2218:329--350, 2001.

\bibitem{chord}
I.~Stoica, R.~Morris, D.~Karger, F.~Kaashoek, and H.~Balakrishnan.
\newblock Chord: {A} scalable {Peer-To-Peer} lookup service for internet
  applications.
\newblock In {\em Proceedings of the 2001 ACM SIGCOMM Conference}, pages
  149--160, 2001.

\end{thebibliography}
\end{document}